\documentclass{amsart}
\usepackage{amsmath,amssymb}
\usepackage{yhmath}

\newcommand{\bdis}{\begin{displaymath}}
\newcommand{\edis}{\end{displaymath}}
\newcommand{\be}{\begin{equation}}
\newcommand{\ee}{\end{equation}}
\newcommand{\mbb}{\mathbb}
\newcommand{\mcal}{\mathcal} 
\newcommand{\mfrak}{\mathfrak}

\newcommand{\vp}{\varphi} 
\newcommand{\vth}{\vartheta}

\newcommand{\zf}{\zeta\left(\frac{1}{2}+it\right)}

\newcommand{\FR}{\frac{x^n+y^n}{z^n}}

\DeclareMathOperator{\im}{Im}

\DeclareMathOperator{\slf}{sl} 
\DeclareMathOperator{\clf}{cl}

\theoremstyle{definition}

\theoremstyle{remark}
\newtheorem{remark}[]{Remark}

\newtheorem*{mydef11}{{\bf Theorem 1}}

\newtheorem*{mydef12}{{\bf Theorem 2}}

\newtheorem*{mydef13}{{\bf Theorem 3}}

\newtheorem*{mydef4}{{\bf Corollary}}

\newtheorem*{mydef51}{{\bf Lemma 1}}

\newtheorem*{mydef52}{{\bf Lemma 2}}

\newtheorem*{mydef81}{{\bf Property 1}}

\newtheorem*{mydef82}{{\bf Property 2}}

\newtheorem*{mydef91}{{\bf Formula 1}}

\numberwithin{equation}{section}

\begin{document}

\title[Jacob's ladders, our old formula (1985)  \dots]{Jacob's ladders, our old formula (1985) and new $\zeta$-equivalent of the Fermat-Wiles theorem on two-parametric set of lemniscates of Bernoulli}

\author{Jan Moser}

\address{Department of Mathematical Analysis and Numerical Mathematics, Comenius University, Mlynska Dolina M105, 842 48 Bratislava, SLOVAKIA}

\email{jan.mozer@fmph.uniba.sk}


\begin{abstract}
In our paper from 1985 we have constructed two integrals of the Riemann's function $Z^2(t)$ over two disconnected sets with asymptotically equal measures such that these two integrals differ by considerably big excess. In the present paper we use the formula for that excess to construct a new $\zeta$-equivalent of the Fermat-Wiles theorem on a two-parametric set of lemniscates of Bernoulli. 
\end{abstract}
\maketitle

\section{Introduction} 

\subsection{} 

Let us remind the Hardy-Littlewood integral\footnote{See \cite{1}} 
\be \label{1.1} 
J=\int_0^T\left|\zf\right|^2{\rm d}t,\ T>0 . 
\ee 
We have obtained in our paper \cite{2} some new facts about the structure of the increments 
\be \label{1.2} 
\int_T^{T+U}\left|\zf\right|^2{\rm d}t,\ U=T^{5/12}\ln^3T 
\ee 
of the Hardy-Littlewood integral as it follows. 

\subsection{} 

At first, let us remind the following notions connected with the Riemann's function 
\be \label{1.3} 
Z(t)=e^{i\vth(t)}\zf;\ Z(t)\in\mbb{R} \ \mbox{if} \ t\in\mbb{R}, 
\ee  
where\footnote{See \cite{11}, (35), (44), (62), comp. \cite{12}, p. 98.}
\be \label{1.4} 
\vth(t)=-\frac t2\ln\pi+\im\ln\Gamma\left(\frac 14+i\frac t2\right)=\frac t2\ln\frac{t}{2\pi}-\frac t2-\frac{\pi}{8}+\mcal{O}\left(\frac 1t\right). 
\ee 
Of course, it is true that 
\be \label{1.5} 
\left|\zf\right|=|Z(t)|. 
\ee 
Let us put (see (\ref{1.4})) 
\bdis 
\vth_1(t)=\frac t2\ln\frac{t}{2\pi}-\frac t2-\frac{\pi}{8}. 
\edis 
In our paper \cite{3} we have defined the following set of sequences 
\be \label{1.6} 
\{g_\nu(\tau)\}_{\nu=1}^\infty,\ \tau\in[-\pi,\pi], 
\ee  
where, see \cite{3}, (6), 
\be \label{1.7} 
\vth_1[g_\nu(\tau)]=\frac{\pi}{2}\nu+\frac{\tau}{2};\ g_\nu(0)=g_\nu. 
\ee 
Next, we have defined\footnote{See \cite{3}, p. 25.} two systems of disconnected sets 
\be \label{1.8} 
\begin{split}
& G_3(v)=G_3(v;T,U)= \\ 
& \bigcup_{T\leq g_{2\nu}\leq T+U}\{t:\ g_{2\nu}(-v)<t<g_{2\nu}(v)\},\ 0\leq v\leq\frac{\pi}{2}, \\ 
& G_4(w)=G_4(w;T,U)= \\
& \bigcup_{T\leq g_{2\nu+1}\leq T+U}\{t:\ g_{2\nu+1}(-w)<t<g_{2\nu+1}(w)\},\ 0\leq w\leq\frac{\pi}{2}. 
\end{split}
\ee 

\subsection{} 

Now, let us put, see (\ref{1.5}), (\ref{1.8}), 
\be \label{1.9} 
\begin{split}
& \sum_{T\leq g_{2\nu}\leq T+U}\int_{g_{2\nu}(-v)}^{g_{2\nu}(v)}Z^2(t){\rm d}t=\int_{G_3(v)}Z^2(t){\rm d}t, \\ 
& \sum_{T\leq g_{2\nu+1}\leq T+U}\int_{g_{2\nu+1}(-w)}^{g_{2\nu+1}(w)}Z^2(t){\rm d}t=\int_{G_4(w)}Z^2(t){\rm d}t. 
\end{split}
\ee 
In this paper we use as a basic formula a variant of our formula \cite{3}, (16). Namely, we use the following result. 

\begin{mydef91}
\be \label{1.10} 
\begin{split}
& \int_{G_3(v)}Z^2(t){\rm d}t-\int_{G_4(w)}Z^2(t){\rm d}t=\\ 
& \frac{4}{\pi}U\sin v	+\mcal{O}(vT^{5/12}\ln^3T)
\end{split}
\ee 
for two continuum systems $\{G_3(v)\}$, $\{G_4(w)\}$ of disconnected sets with almost equal measures\footnote{See \cite{3}, (13).}:
\be \label{1.11} 
m[G_3(v)]=m[G_4(w)]+\mcal{O}\left(\frac{v}{\ln T}\right),\ T\to\infty, 
\ee  
where 
\be \label{1.12} 
U=T^{5/12}\ln^3T,\ 0<v\leq\frac{\pi}{2}. 
\ee 
\end{mydef91} 

\begin{remark}
That is, we have proved in \cite{3} that there is considerable difference between two basic integrals (\ref{1.9}), namely there is a big excess of the first one over the other. 
\end{remark}

\subsection{} 

Now, let us remind the lemniscate of Bernoulli. It is the planar locus traced out by a point $M$ that moves in such a way that the product of its distances from two fixed points 
\be \label{1.13} 
F_1=(a,0),\ F_2=(-a,0)
\ee 
is constant 
\be \label{1.14} 
|F_1M||F_2M|=a^2. 
\ee 

\begin{remark}
Of course, the lemniscate of Bernoulli is a special case of the ovals of Cassini: $a^2\to b^2$ in (\ref{1.14}). 
\end{remark}

Further, let us remind the rationals of Fermat: 
\be \label{1.15} 
\FR,\quad  x,y,z,n\in\mbb{N},\ n\geq 3. 
\ee 

We shall use the two parametric set of the lemniscates of Bernoulli\footnote{See (\ref{1.13}) and (\ref{1.14}).}: 
\be \label{1.16} 
a=\left(\frac{4}{\pi}\sin v\right)^{6/5}\left(\frac{1}{8l_3}\right)^{1/2};\ l_3>0, 0<v\leq \frac{\pi}{2}. 
\ee 

As an example of a new result obtained in this paper we can formulate the following statement: The $\zeta$-condition 
\be \label{1.17} 
\begin{split}
& \lim_{\tau\to\infty}\frac{1}{\tau}
\left\{ 
\sum_{g_{2\nu}\geq\FR\tau}^{g_{2\nu}\leq\FR\tau+U(\FR\tau)}\int_{g_{2\nu}(-v)}^{g_{2\nu}(v)}Z^2(t){\rm d}t- \right. \\ 
& \left. \sum_{g_{2\nu+1}\geq\FR\tau}^{g_{2\nu+1}\leq\FR\tau+U(\FR\tau)}\int_{g_{2\nu+1}(-v)}^{g_{2\nu+1}(v)}Z^2(t){\rm d}t
\right\}^{12/5}\times \\ 
& \left\{ 
\int_{[\FR\tau]^4}^{[\FR\tau+2|F_1M|]^4}\prod_{r=0}^3Z^2(\vp_1^r(t)){\rm d}t\times \right. \\ 
& \left. \int_{[\FR\tau]^3}^{[\FR\tau+2|F_2M|]^3}\prod_{r=0}^2Z^2(\vp_1^r(t)){\rm d}t
\right\}^{-1}\times \\ 
& \left\{ 
\int_{[\FR\tau]^1}^{[\FR\tau+2l_3]^1}Z^2(t){\rm d}t
\right\}^{-1/5}\not=1; \\ 
& [G]^r=\vp_1^{-r}(G), 
\end{split}
\ee 
on the set of all Fermat's rationals and for every fixed 
\be \label{1.18} 
l_3>0,\ 0<v\leq \frac{\pi}{2}, 
\ee 
expresses the $\zeta$-equivalent of the Fermat-Wiles theorem on the two parametric set of corresponding lemniscates of Bernoulli. 

\begin{remark}
Let us remind that C. F. Gauss introduced the lemniscate functions $\slf\tau$ and $\clf\tau$: 
\bdis 
\int_0^{\slf\tau}\frac{{\rm d}x}{\sqrt{1-x^4}}=\tau,\ 
\int_{\clf\tau}^1\frac{{\rm d}x}{\sqrt{1-x^4}}=\tau 
\edis 
as the first elliptic functions on the Gauss plane. 
\end{remark} 

\section{Jacob's ladders: notions and basic geometrical properties}  

\subsection{}

In this paper we use the following notions of our works \cite{5} -- \cite{9}: 
\begin{itemize}
\item[{\tt (a)}] Jacob's ladder $\vp_1(T)$, 
\item[{\tt (b)}] direct iterations of Jacob's ladders 
\bdis 
\begin{split}
	& \vp_1^0(t)=t,\ \vp_1^1(t)=\vp_1(t),\ \vp_1^2(t)=\vp_1(\vp_1(t)),\dots , \\ 
	& \vp_1^k(t)=\vp_1(\vp_1^{k-1}(t))
\end{split}
\edis 
for every fixed natural number $k$, 
\item[{\tt (c)}] reverse iterations of Jacob's ladders 
\be \label{2.1}  
\begin{split}
	& \vp_1^{-1}(T)=\overset{1}{T},\ \vp_1^{-2}(T)=\vp_1^{-1}(\overset{1}{T})=\overset{2}{T},\dots, \\ 
	& \vp_1^{-r}(T)=\vp_1^{-1}(\overset{r-1}{T})=\overset{r}{T},\ r=1,\dots,k, 
\end{split} 
\ee   
where, for example, 
\be \label{2.2} 
\vp_1(\overset{r}{T})=\overset{r-1}{T}
\ee  
for every fixed $k\in\mbb{N}$ and every sufficiently big $T>0$. We also use the properties of the reverse iterations listed below.  
\be \label{2.3}
\overset{r}{T}-\overset{r-1}{T}\sim(1-c)\pi(\overset{r}{T});\ \pi(\overset{r}{T})\sim\frac{\overset{r}{T}}{\ln \overset{r}{T}},\ r=1,\dots,k,\ T\to\infty,  
\ee 
\be \label{2.4} 
\overset{0}{T}=T<\overset{1}{T}(T)<\overset{2}{T}(T)<\dots<\overset{k}{T}(T), 
\ee 
and 
\be \label{2.5} 
T\sim \overset{1}{T}\sim \overset{2}{T}\sim \dots\sim \overset{k}{T},\ T\to\infty.   
\ee  
\end{itemize} 

\begin{remark}
	The asymptotic behaviour of the points 
	\bdis 
	\{T,\overset{1}{T},\dots,\overset{k}{T}\}
	\edis  
	is as follows: at $T\to\infty$ these points recede unboundedly each from other and all together are receding to infinity. Hence, the set of these points behaves at $T\to\infty$ as one-dimensional Friedmann-Hubble expanding Universe. 
\end{remark}  

\subsection{} 

Let us remind that we have proved\footnote{See \cite{9}, (3.4).} the existence of almost linear increments 
\be \label{2.6} 
\begin{split}
& \int_{\overset{r-1}{T}}^{\overset{r}{T}}\left|\zf\right|^2{\rm d}t\sim (1-c)\overset{r-1}{T}, \\ 
& r=1,\dots,k,\ T\to\infty,\ \overset{r}{T}=\overset{r}{T}(T)=\vp_1^{-r}(T)
\end{split} 
\ee 
for the Hardy-Littlewood integral (1918) 
\be \label{2.7} 
J(T)=\int_0^T\left|\zf\right|^2{\rm d}t. 
\ee  

For completeness, we give here some basic geometrical properties related to Jacob's ladders. These are generated by the sequence 
\be \label{2.8} 
T\to \left\{\overset{r}{T}(T)\right\}_{r=1}^k
\ee 
of reverse iterations of the Jacob's ladders for every sufficiently big $T>0$ and every fixed $k\in\mbb{N}$. 

\begin{mydef81}
The sequence (\ref{2.8}) defines a partition of the segment $[T,\overset{k}{T}]$ as follows 
\be \label{2.9} 
|[T,\overset{k}{T}]|=\sum_{r=1}^k|[\overset{r-1}{T},\overset{r}{T}]|
\ee 
on the asymptotically equidistant parts 
\be \label{2.10} 
\begin{split}
& \overset{r}{T}-\overset{r-1}{T}\sim \overset{r+1}{T}-\overset{r}{T}, \\ 
& r=1,\dots,k-1,\ T\to\infty. 
\end{split}
\ee 
\end{mydef81} 

\begin{mydef82}
Simultaneously with the Property 1, the sequence (\ref{2.8}) defines the partition of the integral 
\be \label{2.11} 
\int_T^{\overset{k}{T}}\left|\zf\right|^2{\rm d}t
\ee 
into the parts 
\be \label{2.12} 
\int_T^{\overset{k}{T}}\left|\zf\right|^2{\rm d}t=\sum_{r=1}^k\int_{\overset{r-1}{T}}^{\overset{r}{T}}\left|\zf\right|^2{\rm d}t, 
\ee 
that are asymptotically equal 
\be \label{2.13} 
\int_{\overset{r-1}{T}}^{\overset{r}{T}}\left|\zf\right|^2{\rm d}t\sim \int_{\overset{r}{T}}^{\overset{r+1}{T}}\left|\zf\right|^2{\rm d}t,\ T\to\infty. 
\ee 
\end{mydef82} 

It is clear, that (\ref{2.10}) follows from (\ref{2.3}) and (\ref{2.5}) since 
\be \label{2.14} 
\overset{r}{T}-\overset{r-1}{T}\sim (1-c)\frac{\overset{r}{T}}{\ln \overset{r}{T}}\sim (1-c)\frac{T}{\ln T},\ r=1,\dots,k, 
\ee  
while our eq. (\ref{2.13}) follows from (\ref{2.6}) and (\ref{2.5}).

\section{Some next remarks about the result (1.10) of the paper \cite{3}} 

Since the paper \cite{3} has been published 40 years ago in almost unknown journal, we give here some complements to formulae (\ref{1.8}) -- (\ref{1.12}). 

\subsection{} 

We start with the following remarks. 

\begin{remark}
Since, see (\ref{1.17}), 
\be \label{3.1} 
g_{2\nu}\left(\frac{\pi}{2}\right)=g_{2\nu+1}\left(-\frac{\pi}{2}\right), 
\ee  
then we have the systems of two disconnected sets 
\be \label{3.2} 
\{ G_3(v)\},\ \{ G_4(w)\},\ 0<v,w\leq\frac{\pi}{2} 
\ee  
with the property\footnote{See the symbols of inequalities in curly brackets in (\ref{1.8}).}
\be \label{3.3} 
G_3(v)\cap G_4(w)=\emptyset. 
\ee  
\end{remark} 

\begin{remark}
We have used the symbols $x,y$ in our paper \cite{3} instead of $v,w$ used in the present paper, where $x$ and $y$ are reserved for the contact with the Fermat-Wiles theorem. 
\end{remark} 

\subsection{} 

Since\footnote{See \cite{3}, (11).} 
\be \label{3.4} 
\begin{split}
& g_{2\nu}(v)-g_{2\nu}(-v)=\frac{2v}{\ln\frac{T}{2\pi}}+\mcal{O}\left(\frac{vU}{T\ln^2T}\right), \\ 
& g_{2\nu+1}(w)-g_{2\nu+1}(-w)=\frac{2w}{\ln\frac{T}{2\pi}}+\mcal{O}\left(\frac{wU}{T\ln^2T}\right), 
\end{split}
\ee 
where 
\be\label{3.5} 
g_{2\nu}(-v), g_{2\nu+1}(w)\in [T,T+U],\ 0<v,w\leq \frac{\pi}{2}, 
\ee  
and\footnote{See \cite{2}, (21).} 
\bdis 
\sum_{T\leq g_\nu\leq T+U}1 = \frac{1}{\pi}U\ln\frac{T}{2\pi}+\mcal{O}(1);\ \sum_{g_{2\nu}}=\sum_{g_{2\nu+1}}+O(1), 
\edis 
then for measures of the sets $G_3(v)$ and $G_4(w)$ we have the following formulae 
\be \label{3.6} 
\begin{split}
& m[G_3(v)]=\frac{v}{\pi}U+\mcal{O}\left(\frac{v}{\ln T}\right), \\ 
& m[G_4(v)]=\frac{w}{\pi}U+\mcal{O}\left(\frac{w}{\ln T}\right), 
\end{split}
\ee 
where, of course\footnote{See (\ref{3.5}).},  
\be \label{3.7} 
G_3(v), G_4(w)\subset [T,T+U],\ 0<v,w\leq\frac{\pi}{2}. 
\ee 

\begin{remark}
Let us denote 
\be \label{3.8} 
G_3\left(\frac{\pi}{2}\right)=\bar{G}_3,\ G_4\left(\frac{\pi}{2}\right)=\bar{G}_4. 
\ee 
Then, see (\ref{3.6}), 
\be \label{3.9} 
m[\bar{G}_3]+m[\bar{G}_4]= U+\mcal{O}\left(\frac{1}{\ln T}\right),\ T\to\infty. 
\ee 
\end{remark} 

\subsection{} 

Finally, we present, for completeness, next two results of the paper \cite{3}. Namely, the formulae 
\be \label{3.10} 
\begin{split}
& \int_{G_3(v)}Z^2(t){\rm d}t\sim \int_{G_4(v)}Z^2(t){\rm d}t\sim \frac{v}{\pi}U\ln T,\ T\to\infty, 
\end{split}
\ee   
and also the formula for the difference of the two integrals 
\be \label{3.11} 
\begin{split}
& \int_{G_3(v)}Z^2(t){\rm d}t - \int_{G_4(v)}Z^2(t){\rm d}t\sim \frac{4}{\pi}U\sin v,\ T\to\infty. 
\end{split}
\ee  

\section{The functional generated by the formula (\ref{1.10}) and next $\zeta$-equivalent of the Fermat-Wiles theorem}

\subsection{} 

We rewrite the formula (\ref{1.10}) into the following form, see (\ref{1.5}), (\ref{1.12}), 
\be \label{4.1} 
\begin{split}
& \sum_{g_{2\nu}\geq T}^{g_{2\nu}\leq T+U(T)}\int_{g_{2\nu}(-v)}^{g_{2\nu}(v)}Z^2(t){\rm d}t- \\ 
& \sum_{g_{2\nu+1}\geq T}^{g_{2\nu+1}\leq T+U(T)}\int_{g_{2\nu+1}(-v)}^{g_{2\nu+1}(v)}Z^2(t){\rm d}t = \\ 
& \left\{1+\mcal{O}\left(\frac{1}{\ln T}\right)\right\}\frac{4}{\pi}\sin v T^{5/12}\ln^3T,\ T\to\infty, 
\end{split}
\ee  
and further 
\be \label{4.2} 
\begin{split}
& \left\{
\sum_{g_{2\nu}\geq T}^{g_{2\nu}\leq T+U(T)}\int_{g_{2\nu}(-v)}^{g_{2\nu}(v)}Z^2(t){\rm d}t- \sum_{g_{2\nu+1}\geq T}^{g_{2\nu+1}\leq T+U(T)}\int_{g_{2\nu+1}(-v)}^{g_{2\nu+1}(v)}Z^2(t){\rm d}t
\right\}^{12/5}= \\ 
& \left\{1+\mcal{O}\left(\frac{1}{\ln T}\right)\right\}\left(\frac{4}{\pi}\sin v\right)^{12/5}T\ln^4T\ln^{1/5}T, \\ 
& U(T)=T^{5/12}\ln^3T,\ 0<v\leq\frac{\pi}{2}. 
\end{split}
\ee  
Next, we use in (\ref{4.2}) our formula\footnote{See \cite{10}, (3.18); $f_m\equiv 1$.} 
\be \label{4.3} 
\begin{split}
& 2l\ln^kT=\left\{
\int_{\overset{k}{T}}^{\overset{k}{\wideparen{T+2l}}}\prod_{r=0}^{k-1}Z^2(\vp_1^r(t)){\rm d}t
\right\}\times \left\{
1+\mcal{O}\left(\frac{\ln \ln T}{\ln T}\right)
\right\} , \\ 
& k\in\mbb{N},\ l>0;\ \overset{k}{T}=[T]^k=\vp_1^{-k}(T), \dots
\end{split}
\ee 
for $k=4,3,1$, where we use the partition $7=4+3$ (for example), and obtain the following 
\be \label{4.4} 
\begin{split}
& \left\{
\sum_{g_{2\nu}\geq T}^{g_{2\nu}\leq T+U(T)}\int_{g_{2\nu}(-v)}^{g_{2\nu}(v)}Z^2(t){\rm d}t- \sum_{g_{2\nu+1}\geq T}^{g_{2\nu+1}\leq T+U(T)}\int_{g_{2\nu+1}(-v)}^{g_{2\nu+1}(v)}Z^2(t){\rm d}t
\right\}^{12/5}= \\ 
& \left\{
1+\mcal{O}\left(\frac{\ln \ln T}{\ln T}\right)
\right\}\left(\frac{4}{\pi}\sin v\right)^{12/5}\frac{1}{8l_1l_2l_3}T\times \\ 
& \int_{\overset{4}{T}}^{\overset{4}{\wideparen{T+2l_1}}}\prod_{r=0}^3Z^2(\vp_1^r(t)){\rm d}t\times  \int_{\overset{3}{T}}^{\overset{3}{\wideparen{T+2l_2}}}\prod_{r=0}^2Z^2(\vp_1^r(t)){\rm d}t\times \\ 
& 
\left\{\int_{\overset{1}{T}}^{\overset{1}{\wideparen{T+2l_1}}}Z^2(t){\rm d}t\right\}^{1/5}
\end{split}
\ee 
for every fixed 
\be \label{4.5} 
l_1,l_2,l_3>0, 0<v\leq\frac{\pi}{2}. 
\ee 

\subsection{} 

Next, we introduce the following condition, see (\ref{4.4}), 
\be \label{4.6} 
\left(\frac{4}{\pi}\sin v\right)^{12/5}\frac{1}{8l_1l_2l_3}=1, 
\ee 
that is, for example, 
\be \label{4.7} 
l_1l_2=\left(\frac{4}{\pi}\sin v\right)^{12/5}\frac{1}{8l_3},\ l_1,l_2,l_3>0,\ 0<v\leq \frac{\pi}{2}, 
\ee  
or 
\be \label{4.8} 
l_1l_2=a^2,\ a=\left(\frac{4}{\pi}\sin v\right)^{6/5}\left(\frac{1}{8l_3}\right)^{1/2}. 
\ee 
Now, we put, for example\footnote{See (\ref{1.13}), (\ref{1.14}), (\ref{1.16}).} 
\be \label{4.9} 
l_1=|F_1M|,\ l_2=|F_2M| 
\ee 
into (\ref{4.8}) and we obtain\footnote{See (\ref{4.4}), (\ref{4.6}) -- (\ref{4.8}).} the following result. 

\begin{mydef51}
It is true that 
\be \label{4.10} 
\begin{split}
& \left\{
\sum_{g_{2\nu}\geq T}^{g_{2\nu}\leq T+U(T)}\int_{g_{2\nu}(-v)}^{g_{2\nu}(v)}Z^2(t){\rm d}t- \sum_{g_{2\nu+1}\geq T}^{g_{2\nu+1}\leq T+U(T)}\int_{g_{2\nu+1}(-v)}^{g_{2\nu+1}(v)}Z^2(t){\rm d}t
\right\}^{12/5}=\\ 
& \left\{
1+\mcal{O}\left(\frac{\ln \ln T}{\ln T}\right)
\right\}T\int_{\overset{4}{T}}^{\overset{4}{\wideparen{T+2|F_1M|}}}\prod_{r=0}^3Z^2(\vp_1^r(t)){\rm d}t\times \\ 
& \int_{\overset{3}{T}}^{\overset{3}{\wideparen{T+2|F_2M|}}}\prod_{r=0}^2Z^2(\vp_1^r(t)){\rm d}t\times 
\left\{
\int_{\overset{1}{T}}^{\overset{1}{\wideparen{T+2l_3}}}Z^2(t){\rm d}t
\right\}^{1/5}
\end{split}
\ee 
for every fixed 
\be \label{4.11} 
l_3>0,\ 0<v\leq\frac{\pi}{2}, 
\ee  
and every point $M\in\mcal{L}[l_3,v]$ i.e. arbitrary point of the lemniscate of Bernoulli (\ref{4.8}) $=\mcal{L}[l_3,v]$.
\end{mydef51}

\subsection{} 

Next, we make the substitution 
\be \label{4.12} 
T=x\tau,\ x>0;\ \{T\to+\infty\} \ \Leftrightarrow \ \{\tau\to+\infty\}
\ee 
in (\ref{4.10}) and we obtain, as usually in our theory, the following functional. 

\begin{mydef11}
It is true that 
\be \label{4.13} 
\begin{split}
& \lim_{\tau\to\infty}\frac{1}{\tau}\times \\ 
& 
\left\{
\sum_{g_{2\nu}\geq x\tau}^{g_{2\nu}\leq x\tau+U(x\tau)}\int_{g_{2\nu}(-v)}^{g_{2\nu}(v)}Z^2(t){\rm d}t- \sum_{g_{2\nu+1}\geq x\tau}^{g_{2\nu+1}\leq x\tau+U(x\tau)}\int_{g_{2\nu+1}(-v)}^{g_{2\nu+1}(v)}Z^2(t){\rm d}t
\right\}^{12/5}\times \\ 
& \left\{
\int_{[x\tau]^4}^{[x\tau+2|F_1M|]^4}\prod_{r=0}^3Z^2(\vp_1^r(t)){\rm d}t\times 
\int_{[x\tau]^3}^{[x\tau+2|F_2M|]^3}\prod_{r=0}^2Z^2(\vp_1^r(t)){\rm d}t
\right\}^{-1}\times \\ 
& \left\{\int_{[x\tau]^1}^{[x\tau+2l_3]^1}Z^2(t){\rm d}t\right\}^{-1/5}=x 
\end{split}
\ee 
for every fixed 
\be \label{4.14} 
x>0,\ l_3>0,\ 0<v\leq \frac{\pi}{2}, 
\ee 
and for every point $M\in\mcal{L}[l_3,v]$. 
\end{mydef11} 

\subsection{} 

We obtain from (\ref{4.13}) the following consequence for Fermat's rationals\footnote{See (\ref{1.15}).} 

\begin{mydef4}
\be \label{4.15} 
\begin{split}
	& \lim_{\tau\to\infty}\frac{1}{\tau}\times \\ 
	& 
	\left\{
	\sum_{g_{2\nu}\geq \FR\tau}^{g_{2\nu}\leq \FR\tau+U(\FR\tau)}\int_{g_{2\nu}(-v)}^{g_{2\nu}(v)}Z^2(t){\rm d}t- \right. \\  
	& \left. \sum_{g_{2\nu+1}\geq \FR\tau}^{g_{2\nu+1}\leq \FR\tau+U(\FR\tau)}\int_{g_{2\nu+1}(-v)}^{g_{2\nu+1}(v)}Z^2(t){\rm d}t
	\right\}^{12/5}\times \\ 
	& \left\{
	\int_{[\FR\tau]^4}^{[\FR\tau+2|F_1M|]^4}\prod_{r=0}^3Z^2(\vp_1^r(t)){\rm d}t\times \right. \\  
	& \left. \int_{[\FR\tau]^3}^{[\FR\tau+2|F_2M|]^3}\prod_{r=0}^2Z^2(\vp_1^r(t)){\rm d}t
	\right\}^{-1}\times \\ 
	& \left\{\int_{[\FR\tau]^1}^{[\FR\tau+2l_3]^1}Z^2(t){\rm d}t\right\}^{-1/5}=\FR
\end{split}
\ee 
for every fixed 
\be \label{4.16} 
l_3>0,\ 0<v\leq\frac{\pi}{2}, 
\ee  
and for every point $M\in\mcal{L}[l_3,v]$. 
\end{mydef4}

Consequently, we obtain from (\ref{4.15}) the following result. 

\begin{mydef12}
The $\zeta$-condition 
\be \label{4.17} 
\begin{split}
& \lim_{\tau\to\infty}\frac{1}{\tau}\times \\ 
& 
\left\{
\sum_{g_{2\nu}\geq \FR\tau}^{g_{2\nu}\leq \FR\tau+U(\FR\tau)}\int_{g_{2\nu}(-v)}^{g_{2\nu}(v)}Z^2(t){\rm d}t- \right. \\  
& \left. \sum_{g_{2\nu+1}\geq \FR\tau}^{g_{2\nu+1}\leq \FR\tau+U(\FR\tau)}\int_{g_{2\nu+1}(-v)}^{g_{2\nu+1}(v)}Z^2(t){\rm d}t
\right\}^{12/5}\times \\ 
& \left\{
\int_{[\FR\tau]^4}^{[\FR\tau+2|F_1M|]^4}\prod_{r=0}^3Z^2(\vp_1^r(t)){\rm d}t\times \right. \\  
& \left. \int_{[\FR\tau]^3}^{[\FR\tau+2|F_2M|]^3}\prod_{r=0}^2Z^2(\vp_1^r(t)){\rm d}t
\right\}^{-1}\times \\ 
& \left\{\int_{[\FR\tau]^1}^{[\FR\tau+2l_3]^1}Z^2(t){\rm d}t\right\}^{-1/5}\not=1
\end{split}
\ee 
on the set of all Fermat's rationals and for every fixed 
\be \label{4.18} 
l_3>0,\ 0<v\leq \frac{\pi}{2}, 
\ee  
and also for every point $M\in\mcal{L}[l_3,v]$, expresses the new $\zeta$-equivalent of the Fermat-Wiles theorem (on two parametric set of lemniscates of Bernloulli, see (\ref{1.13}), (\ref{1.14}), (\ref{1.16}).)
\end{mydef12}

\section{A factorization formula for the excess (\ref{1.10})} 

\subsection{} 

Next, we shall use in (\ref{4.4}) our almost linear formula, see \cite{8}, (3.4), (3.6), $r=1$, in the form 
\be \label{5.1} 
(1-c)T=\{1+\mcal{O}(T^{-1/3+\delta})\}\int_T^{\overset{1}{T}}Z^2(t){\rm d}t, 
\ee  
that gives us the following formula. 

\begin{mydef52}
\be \label{5.2} 
\begin{split}
& \left\{
\sum_{g_{2\nu}\geq T}^{g_{2\nu}\leq T+U(T)}\int_{g_{2\nu}(-v)}^{g_{2\nu}(v)}Z^2(t){\rm d}t- \sum_{g_{2\nu+1}\geq T}^{g_{2\nu+1}\leq T+U(T)}\int_{g_{2\nu+1}(-v)}^{g_{2\nu+1}(v)}Z^2(t){\rm d}t
\right\}^{12/5}= \\ 
& \left\{
1+\mcal{O}\left(\frac{\ln \ln T}{\ln T}\right)
\right\}\left(\frac{4}{\pi}\sin v\right)^{12/5}\frac{1}{8(1-c)l_1l_2l_3}\times \int_T^{\overset{1}{T}}Z^2(t){\rm d}t\times \\ 
& \left\{\int_{\overset{1}{T}}^{\overset{1}{\wideparen{T+2l_3}}}Z^2(t){\rm d}t\right\}^{1/5}\times \int_{\overset{3}{T}}^{\overset{3}{\wideparen{T+2l_2}}}\prod_{r=0}^2Z^2(\vp_1^r(t)){\rm d}t \times \\  
& \int_{\overset{4}{T}}^{\overset{4}{\wideparen{T+2l_3}}}\prod_{r=0}^3Z^2(\vp_1^r(t)){\rm d}t
\end{split}
\ee 
for every fixed 
\be \label{5.3} 
l_1,l_2,l_3>0,\ 0<v\leq \frac{\pi}{2}. 
\ee 
\end{mydef52}

\subsection{} 

Now, we use the following condition for the asymptotic equality\footnote{Comp. (\ref{4.6}).} 
\be \label{5.4} 
\left(\frac{4}{\pi}\sin v\right)^{12/5}\frac{1}{8(1-c)l_1l_2l_3}=1
\ee 
i.e., for example\footnote{Comp. (\ref{4.7}), (\ref{4.8}).}, 
\be \label{5.5} 
l_1l_2\bar{a}^2,\ \bar{a}=
\left(\frac{4}{\pi}\sin v\right)^{6/5}\left(\frac{1}{8(1-c)l_3}\right)^{1/2}. 
\ee  
Consequently, we put, for example\footnote{Comp. (\ref{1.13}) -- (\ref{1.16}).} 
\be \label{5.6} 
F_3=F_3(\bar{a},0),\ F_4=F_4(-\bar{a},0)
\ee  
and after this we substitute the values 
\be \label{5.7} 
l_1=|F_3M|,\ l_2=|F_4M|
\ee 
into (\ref{5.2}) and we obtain\footnote{Comp. (\ref{4.6}) -- (\ref{4.9}).} the following result. 

\begin{mydef13}
\be \label{5.8} 
\begin{split}
& \sum_{g_{2\nu}\geq T}^{g_{2\nu}\leq T+U(T)}\int_{g_{2\nu}(-v)}^{g_{2\nu}(v)}Z^2(t){\rm d}t- \sum_{g_{2\nu+1}\geq T}^{g_{2\nu+1}\leq T+U(T)}\int_{g_{2\nu+1}(-v)}^{g_{2\nu+1}(v)}Z^2(t){\rm d}t\sim \\ 
& \left\{\int_T^{\overset{1}{T}}Z^2(t){\rm d}t\right\}^{5/12}\times 
\left\{
\int_{\overset{1}{T}}^{\overset{1}{\wideparen{T+2l_3}}}Z^2(t){\rm d}t
\right\}^{1/12}\times \\ 
& \left\{
\int_{\overset{3}{T}}^{\overset{3}{\wideparen{T+2|F_4M|}}}\prod_{r=0}^2Z^2(\vp_1^r(t)){\rm d}t
\right\}^{5/12} \times \\ 
& \left\{
\int_{\overset{4}{T}}^{\overset{4}{\wideparen{T+2|F_3M|}}}\prod_{r=0}^3Z^2(\vp_1^r(t)){\rm d}t
\right\}^{5/12}
\end{split}
\ee 
for every fixed 
\be \label{5.9} 
l_3>0,\ 0<v\leq \frac{\pi}{2}, 
\ee  
and every point $M\in\bar{\mcal{L}}[l_3,v]$, i.e. the corresponding lemniscate of Bernoulli (5.5) = $\bar{\mcal{L}}[l_3,v]$. 
\end{mydef13} 

\begin{remark}
The formula (\ref{5.8}) expresses  the asymptotic factorization  of the excess (\ref{1.10}) by means of the basic set of integrals 
\be \label{5.10} 
\begin{split}
& \left\{ 
\int_T^{[T]^1}Z^2(t){\rm d}t,\ \int_{[T]^1}^{[T+2l_3]^1}Z^2(t){\rm d}t \right. \\ 
& \left. \int_{[T]^3}^{[T+2l_2]^3}\prod_{r=0}^2 Z^2(\vp_1^r(t)){\rm d}t, \ 
\int_{[T]^4}^{[T+2l_1]^4}\prod_{r=0}^3 Z^2(\vp_1^r(t)){\rm d}t
\right\}, 
\end{split}
\ee 
where, of course, 
\bdis 
[G]^r=\vp_1^r(G), 
\edis  
on the other two-parametric system of lemniscates of Bernoulli $\bar{\mcal{L}}[l_3,v]$. 
\end{remark}

\section{On the structure of the set of intervals under question} 

\subsection{} 

Our results listed above contain integrations over the intervals 
\be \label{6.1} 
\begin{split}
& \{(T,[T]^1),\ ([T]^1,[T+2l_3]^1),\ ([T]^3,[T+2l_2]^3), \\ 
& ([T]^4,[T+2l_1]^4)\}. 
\end{split}
\ee  
We use the formula\footnote{Comp. (\ref{2.14}).}  
\be \label{6.2} 
[T]^r-[T]^{r-1}\sim (1-c)\frac{T}{\ln T},\ r=1,\dots,k
\ee  
where 
\be \label{6.3} 
[T]^r=\overset{r}{T}=\vp_1^{-r}(T), 
\ee  
in the following cases: 
\begin{itemize}
	\item[(a)] Since the summation of (\ref{6.2}) for $r=1,2,3$ gives 
	\be \label{6.4} 
	[T]^3-T\sim 3(1-c)\frac{T}{\ln T}, 
	\ee  
	and 
	\be \label{6.5} 
	[T]^3-[T+2l_3]^1-T-2l_3\sim (1-c)\frac{T}{\ln T}, 
	\ee 
	then we obtain  the following difference 
	\be \label{6.6} 
	[T]^3-[T+2l_3]^1-T-2l_3\sim 2(1-c)\frac{T}{\ln T},\ T\to\infty, 
	\ee  
	for every fixed $l_3>0$. Consequently, it is true that 
	\be \label{6.7} 
	[T]^3>[T+2l_3]^1,\ T\to\infty, 
	\ee 
	and the inequality holds true with a big reserve. 
	\item[(b)] Since the summation of (\ref{6.2}) for $r=1,2,3,4$ gives 
	\be \label{6.8} 
	[T]^4-T\sim 4(1-c)\frac{T}{\ln T}, 
	\ee  
	and the summation for $r=1,2,3$ with the substitution 
	\be \label{6.9} 
	T\to T+2l_2;\ \frac{T+2l_2}{\ln(T+2l_2)}\sim \frac{T}{\ln T},\ T\to\infty, 
	\ee  
	for every fixed $l_2>0$ gives 
	\be \label{6.10} 
	[T+2l_2]^3-T\sim 2l_2+3(1-c)\frac{T}{\ln T}\sim 3(1-c)\frac{T}{\ln T}, 
	\ee  
	then we obtain\footnote{See (\ref{6.8}), (\ref{6.10}).} 
	\be \label{6.11} 
	[T]^4-[T+2l_2]^3\sim (1-c)\frac{T}{\ln T}. 
	\ee  
	Consequently, it is true that 
	\be \label{6.12} 
	[T]^4>[T+2l_2]^3,\ T\to\infty, 
	\ee 
	and the inequality holds true with a big reserve. 
	\item[(c)] Only $0$-adjacent interval corresponds to the first and the second of intervals in (\ref{6.1}). 
\end{itemize} 

\begin{remark}
It is true that the set of intervals (\ref{6.1}) is the disconnected set  with the big measures of the second and the third adjacent intervals, comp. (\ref{6.6}) and (\ref{6.11}).  
\end{remark} 

\section{Concluding remarks} 

\subsection{} 

Let us remind that the oval of Cassini is the locus traced out by a point $M$ moving so that the product of its distances from two given points $F_1$ and $F_2$ is constant 
\be \label{7.1} 
|F_1M||F_2M|=b^2, 
\ee  
where 
\be \label{7.2} 
F_1=F_1(a,0),\ F_2=F_2(-a,0). 
\ee 
\begin{itemize}
	\item[($\alpha$)] In the critical case $a=b$ the curve is known as lemniscate of Bernoulli. 
	\item[($\beta$)] In the case $b<a$ the curve consists of two separated ovals. Each of points $F_1$, $F_2$ is surrounded by one of these ovals. 
	 \item[($\gamma$)] In the case $b>a$ the curve is an oval with both pints $F_1,F_2$ in its interior. 
\end{itemize} 

\subsection{} 

Since we wish to keep $a^2$ with 
\be \label{7.3} 
a=\left(\frac{4}{\pi}\sin v\right)^{6/5}\left(\frac{1}{8l_3}\right)^{1/2},\ l_3>0,\ 0<v\leq \frac{\pi}{2}
\ee  
on the right-hand side of (\ref{7.1}), we put 
\be \label{7.4} 
|F_5M||F_6M|=a^2, 
\ee  
where 
\be \label{7.5} 
F_5=F_5(c,0),\ F_6=F_6(-c,0). 
\ee 
We, however, do not with to introduce a new parameter $c$, so we shall consider only following two possibilities 
\be \label{7.6} 
c=\frac 12a,\ c=2a 
\ee 
in (\ref{7.4}), (\ref{7.5}). 

\begin{remark}
Since\footnote{See (\ref{7.3}).} 
\be \label{7.7} 
a>0, 
\ee  
then we have 
\begin{itemize}
	\item[(a)] for 
	\be \label{7.8} 
	c=\frac 12 a \ \to   \ \mfrak{O}_1[l_3,v]
	\ee 
	the continuum set of Cassini ovals $\mfrak{O}_1[l_3,v]$ of the type ($\gamma$), 
	\item[(b)] for 
	\be \label{7.9} 
	c= 2a \ \to \ \mfrak{O}_2[l_3,v]
	\ee 
	the continuum set of Cassini ovals $\mfrak{O}_2[l_3,v]$ of the type ($\beta$). 
\end{itemize}
\end{remark} 

\begin{remark}
Finally, if we substitute the values 
\be \label{7.10} 
|F_7M|,\ |F_8M|
\ee 
with 
\be \label{7.11} 
F_7=F_7(\frac 12 a,0),\ F_8=F_8(-\frac 12 a,0), 
\ee  
and also substitute the values 
\be \label{7.12} 
|F_9M|,\ |F_{10}M|
\ee 
with 
\be \label{7.13} 
F_9=F_9(2a,0),\ F_{10}=F_{10}(-2a,0), 
\ee 
in (\ref{4.13}) and (\ref{4.17}) instead of $|F_1M|$ and $|F_2M|$, respectively, and next in (\ref{5.8}) instead of $|F_3M|$ and $|F_4M|$, then we obtain the corresponding new Theorem on the continuum sets of Cassini ovals $\mfrak{O}_1[l_3,v]$ and $\mfrak{O}_2[l_3,v]$. 
\end{remark}

I would like to thank Michal Demetrian for his moral support of my study of Jacob's ladders.

\end{document}